\documentclass{proc-l}
\usepackage{amssymb,float,graphicx}

\newtheorem{theorem}{Theorem}[section]
\newtheorem{lemma}[theorem]{Lemma}
\newtheorem{proposition}[theorem]{Proposition}
\newtheorem{corollary}[theorem]{Corollary}
\theoremstyle{remark}
\newtheorem*{remark}{Remark}
\numberwithin{equation}{section}

\newcommand{\R}{\mathbb{R}}
\newcommand{\dd}{\,\mathrm{d}}
\newcommand{\one}{\mathbf{1}}
\newcommand{\Sch}{\mathcal{S}}
\newcommand{\F}{\mathcal{F}}
\DeclareMathOperator{\re}{Re}
\DeclareMathOperator{\im}{Im}
\DeclareMathOperator{\li}{li}

\title[Non-divisibility of modified zeta distributions]{Failure of quasi-infinite divisibility for modified Riemann zeta distributions}
\author{Oscar L\'opez-Alfonso}
\address{Departamento de Matem\'aticas, Facultad de Ciencias,
\newline Universidad Nacional de Colombia, Sede Bogot\'a, Colombia}
\email{ojlopeza@unal.edu.co}
\thanks{This work was supported by Universidad Nacional de Colombia through the author's sabbatical leave (Resoluci\'on 2117 of 2026, Facultad de Ciencias, Sede Bogot\'a).}
\copyrightinfo{}{}
\makeatletter
\AtBeginDocument{\let\@serieslogo\@empty}
\makeatother
\begin{document}

\begin{abstract}
We prove that, for every $0<\theta<1$, the normalized Riemann zeta function $t\mapsto\zeta(\theta-it)/\zeta(\theta)$ is not a quotient of infinitely divisible characteristic functions. Consequently no quasi-infinitely divisible characteristic factor makes Nakamura's modified zeta law quasi-infinitely divisible. The proof combines the L\'evy--Khintchine formula, Voronin's universality theorem, and the functional equation. We also give an obstruction for characteristic factors with a finite exponential moment at one point. If $\theta$ exceeds the supremum $\Theta$ of the real parts of the nontrivial zeros of $\zeta$, we derive an improper L\'evy--Khintchine-type representation with atoms at logarithms of prime powers and explicit continuous terms. This hypothesis is sharp except possibly at $\theta=\Theta$: if the improper integral converges at a single $t\ne0$, then $\theta\ge\tfrac12$ and $\zeta(s)\ne0$ for $\operatorname{Re}s>\theta$.
\end{abstract}

\maketitle

\section{Introduction}

For $\theta>1$, the Euler product shows that $\zeta(\theta-it)/\zeta(\theta)$ is a compound Poisson characteristic function. For $0<\theta<1$, the alternating series for $(1-2^{1-\theta})\zeta(\theta)$ gives $\zeta(\theta)<0$; in particular $\zeta(\theta)\ne0$. In this strip set
\begin{equation}\label{eq:notation}
Z_\theta(t):=\frac{\zeta(\theta-it)}{\zeta(\theta)},
\qquad
F_\theta(t):=Z_\theta(t)\frac{\theta}{\theta-it}.
\end{equation}
Nakamura \cite[Theorem~1.2]{Nak15a} proved that $F_\theta$ is a characteristic function, and he calls its law the modified Riemann zeta distribution. The factor $\theta/(\theta-it)$, the characteristic function of an exponential law with rate $\theta$, compensates for the growth of $\zeta(\theta-it)$ in $t$. More generally, for real $c\ge1$ let
\[
\varphi_{\theta,c}(t):=Z_\theta(t)\Bigl(\frac{\theta}{\theta-it}\Bigr)^{c}=F_\theta(t)\Bigl(\frac{\theta}{\theta-it}\Bigr)^{c-1},
\]
where the power is defined by the principal logarithm, which is continuous because $\re(\theta-it)>0$. For $c>1$, the last factor is the characteristic function of a gamma law with shape $c-1$ and rate $\theta$, so $\varphi_{\theta,c}$ is the characteristic function of $X+Y$, where $X$ has characteristic function $F_\theta$ and $Y$ is an independent gamma variable. For $c=1$, $\varphi_{\theta,1}=F_\theta$.

A characteristic function is \emph{quasi-infinitely divisible} (QID) if it is a quotient of two infinitely divisible (ID) characteristic functions \cite[Remark~2.4]{LPS18}. Such quotients form a group under multiplication. Nakamura \cite[Theorem~3.2]{Nak15a} showed that the Riemann hypothesis holds if $F_\theta$ is ID for every $\tfrac12<\theta<1$. In \cite[Remark~3.3]{Nak15a} he drew on the value distribution of $\zeta$ to question the converse. He suggested that, under the Riemann hypothesis, $\varphi_{\theta,c}$ might be ID for large integers $c$, or that $GF_\theta$ might be ID for a suitable ID characteristic function $G$. Feuerverger and Martin \cite[\S~10.4, item~(7), pp.~270--271]{FM26} recall the gamma-factor question.

\begin{theorem}\label{thm:main}
For every $0<\theta<1$, the function $Z_\theta$ is not a quotient of ID characteristic functions. In particular, for any QID characteristic function $G$, neither $GZ_\theta$ nor $GF_\theta$ is a QID characteristic function, and $\varphi_{\theta,c}$ is not QID for any $c\ge1$.
\end{theorem}

Since ID characteristic functions are QID, Theorem~\ref{thm:main} shows that neither suggestion in \cite[Remark~3.3]{Nak15a} can hold, with or without the Riemann hypothesis, and that the hypothesis of \cite[Theorem~3.2]{Nak15a} is never satisfied.

Theorem~\ref{thm:main} is proved in Section~\ref{sec:main}, where we also deduce that Nakamura's complete Riemann zeta distribution is not QID for any $0<\theta<1$ (Corollary~\ref{cor:xi}). The proof of Theorem~\ref{thm:main} rests on the L\'evy--Khintchine formula, Voronin's universality theorem, and the functional equation. A smoothed logarithm of a quotient of ID characteristic functions is bounded (Lemma~\ref{lem:smoothing}), whereas universality forces $\log|Z_\theta(t)|$ to follow arbitrarily steep parabolas on short intervals. The functional equation reduces the range $0<\theta<\tfrac12$ to $\tfrac12<\theta<1$, and a zero on the critical line settles $\theta=\tfrac12$.

The remaining sections complement Theorem~\ref{thm:main}. In Section~\ref{sec:tilt} the QID factor is replaced by an arbitrary characteristic function $H$ with a finite exponential moment at one point, and $HF_\theta$ is shown not to be ID (Theorem~\ref{thm:tilt}). Section~\ref{sec:primes} explains the obstruction arithmetically. If $\theta>\Theta$, where $\Theta$ is the supremum of the real parts of the nontrivial zeros of $\zeta$, the logarithm of $Z_\theta$ has an improper L\'evy--Khintchine-type representation whose signed measure consists of atoms at the logarithms of prime powers and a negative density that nothing compensates on the intervals $(\log m,\log(m+1))$ (Theorem~\ref{thm:primes}). This yields a second proof of Theorem~\ref{thm:main} for such $\theta$. Proposition~\ref{prop:sharp} shows that the representation fails at every $t\ne0$ when $\theta<\Theta$, so the hypothesis is sharp except possibly at $\theta=\Theta$. Section~\ref{sec:numerics} checks the representation numerically, using only values of $\zeta$.

\section{Proof of the quotient obstruction}\label{sec:main}

For a locally integrable function $f$ and $h>0$, put
\[
 D_h f(t):=f(t)-\frac{1}{2h}\int_{-h}^h f(t+u)\dd u,
 \qquad w_h(x):=1-\frac{\sin(hx)}{hx},
 \qquad w_h(0):=0.
\]

\begin{lemma}\label{lem:smoothing}
If $\chi$ is an ID characteristic function with Gaussian
coefficient $a$ and L\'evy measure $\nu$, then its distinguished logarithm
(the continuous logarithm with value zero at $t=0$) satisfies
\[
 D_h\log\chi(t)=\frac{ah^2}{6}+
                  \int_{\R}e^{itx}w_h(x)\,\nu(\mathrm{d}x).
\]
If $Q=\chi_1/\chi_2$ with $\chi_1,\chi_2$ ID, then $D_h\log Q$ is the
Fourier transform of a finite signed measure in the $e^{itx}$ convention.
In particular, $D_h\log|Q|$ is bounded.
\end{lemma}

\begin{proof}
Use the L\'evy--Khintchine formula with Gaussian term $-at^2/2$ and
truncation $x\one_{\{|x|\le1\}}$. The averaging calculation is as in
\cite[Theorem~8.1 and (8.10)]{Sat99}. The operator $D_h$ annihilates constants and linear functions, satisfies
$D_h(t^2)=-h^2/3$, and sends $e^{itx}$ to $e^{itx}w_h(x)$. In the
L\'evy--Khintchine formula the integrand at $t+u$, for $|u|\le h$, is
bounded in absolute value by a constant depending on $t$ and $h$ times
$1\wedge x^2$. Fubini's theorem therefore gives the displayed identity.
The bound $0\le w_h(x)\le\min\{2,h^2x^2/6\}$ makes $w_h\nu$ finite.
By normalization at zero, $\log Q=\log\chi_1-\log\chi_2$.
Subtracting the two Fourier transforms gives the claim for $Q$.
\end{proof}

\begin{proof}[Proof of Theorem~\ref{thm:main}]
First suppose that $\tfrac12<\theta<1$. Choose $r$ with $|\theta-\tfrac34|<r<\tfrac14$, and then $h>0$ so small that the segment
$\{\,\theta+iu:|u|\le h\,\}$ lies in the disk
$|s-\tfrac34|\le r$. For $L>0$, the entire function
\[
 f_L(s):=\exp\!\left(-L\frac{(s-\theta)^2}{h^2}\right)
\]
has no zeros. By Voronin's universality theorem \cite{Vor75}, there is a
vertical shift $\tau$ for which $|\zeta(s+i\tau)-f_L(s)|<1/2$ throughout
the disk. On the segment we have $|f_L(\theta+iu)|=\exp(Lu^2/h^2)\ge1$.
Consequently,
\[
 \left|\log|\zeta(\theta+i(\tau+u))|-\frac{Lu^2}{h^2}\right|
 \le \log 2,\qquad |u|\le h.
\]
If $Z_\theta$ were a quotient of the required form, Lemma~\ref{lem:smoothing}
would make $D_h\log|Z_\theta|$ bounded. However, by conjugation,
$\log|Z_\theta(t)|=\log|\zeta(\theta+it)|-\log|\zeta(\theta)|$.
Averaging the last estimate therefore gives
\[
 \left|D_h\log|Z_\theta|(\tau)+\frac L3\right|\le 2\log 2.
\]
Since $L$ is arbitrary, this contradicts boundedness.

For $0<\theta<\tfrac12$, the functional equation \cite[Chap.~II]{Tit86} gives
\[
 Z_\theta(t)=\pi^{-it}
   \frac{\Gamma((1-\theta+it)/2)}{\Gamma((1-\theta)/2)}
   \frac{\Gamma(\theta/2)}{\Gamma((\theta-it)/2)}
   Z_{1-\theta}(-t).
\]
The factor $\pi^{-it}$ is the characteristic function of a point mass.
For $\alpha>0$, the functions $\Gamma(\alpha\pm it/2)/\Gamma(\alpha)$ are ID
characteristic functions, by the infinite divisibility of
the log-gamma law \cite[Lemma~2.8]{Nak15b} and invariance under scaling
and reflection. Quotients of ID characteristic functions
are closed under multiplication, reciprocation, and reflection. Hence
such a representation for $Z_\theta$ would give one for $Z_{1-\theta}$,
which has already been ruled out. For $\theta=\tfrac12$, a zero on the
critical line \cite{Har14} gives a zero of $Z_\theta$. A quotient
of ID characteristic functions has no zeros.

If $GZ_\theta$ were QID, division by $G$ would give the excluded
representation. The other claims follow because the factors relating
$GF_\theta$ and $\varphi_{\theta,c}$ to $Z_\theta$ are products of QID
characteristic functions and ID gamma characteristic functions.
\end{proof}

Theorem~\ref{thm:main} also applies to two other zeta distributions of Nakamura.

\begin{remark}
Nakamura \cite[Theorem~2.1(1) with $a=z=1$]{Nak16} showed that
$\Gamma(\theta+it)\zeta(\theta+it)/(\Gamma(\theta)\zeta(\theta))$ is a
characteristic function for $0<\theta<1$. Its gamma factor is ID by
\cite[Lemma~2.8]{Nak15b}. If the law were QID, division by that factor
and reflection would express $Z_\theta$ as a quotient of ID characteristic
functions, contrary to Theorem~\ref{thm:main}.
\end{remark}

\begin{corollary}\label{cor:xi}
For every $0<\theta<1$, the characteristic function $\Xi_\theta(t):=\xi(\theta-it)/\xi(\theta)$ of Nakamura's complete Riemann zeta distribution \cite[Theorem~1.1]{Nak15b} is not QID.
\end{corollary}

\begin{proof}
With $\xi(s)=\tfrac12s(s-1)\pi^{-s/2}\Gamma(s/2)\zeta(s)$ we have
\[
\Xi_\theta(t)=\pi^{it/2}\,\frac{\Gamma((\theta-it)/2)}{\Gamma(\theta/2)}\,
\frac{\theta-it}{\theta}\,\frac{1-\theta+it}{1-\theta}\,Z_\theta(t).
\]
The first two factors are ID characteristic functions: the first is that of a point mass, and the second is a rescaled log-gamma characteristic function \cite[Lemma~2.8]{Nak15b}. The next two are the reciprocals of the characteristic functions $\theta/(\theta-it)$ and $(1-\theta)/(1-\theta+it)$ of exponential laws on $(0,\infty)$ and $(-\infty,0)$. If $\Xi_\theta$ were QID, then $Z_\theta$ would be a quotient of ID characteristic functions, contrary to Theorem~\ref{thm:main}.
\end{proof}

Nakamura \cite[Theorem~1.2]{Nak15b} showed that $\Xi_\theta$ is pretended-infinitely divisible for every $\tfrac12<\theta<1$ if and only if the Riemann hypothesis holds: its logarithm has the L\'evy--Khintchine form with a signed measure, but without the integrability condition $\int(|x|^2\wedge1)\,|\nu|(\mathrm{d}x)<\infty$. By Corollary~\ref{cor:xi}, for no $0<\theta<1$ does $\Xi_\theta$ admit such a representation that also satisfies this condition.

\section{An exponential-moment obstruction}\label{sec:tilt}

If a law has finite exponential moments on an interval, its characteristic function extends analytically to the corresponding strip, and when the law is ID the extension has no zeros there \cite[Theorem~4]{LS52}. The following result requires a finite exponential moment at a single point $\lambda$ only. An exponential tilt by $e^{\lambda y}$ carries the line $\re z=\lambda$ of the moment generating function $z\mapsto\int_\R e^{zy}\,\mu(\mathrm{d}y)$ to the real axis and preserves infinite divisibility, so it suffices to produce a real zero of the characteristic function of the tilted law.

\begin{theorem}\label{thm:tilt}
Let $0<\theta<1$, let $\rho$ be a nontrivial zero of $\zeta$, and put
$\lambda=\theta-\re\rho$. If $H$ is the characteristic function of a law
$\eta$ with $\int_\R e^{\lambda y}\eta(\mathrm{d}y)<\infty$, then $HF_\theta$ is
the characteristic function of a law that is not ID.
\end{theorem}

\begin{proof}
Write $M_\mu(\lambda)=\int_\R e^{\lambda y}\mu(\mathrm{d}y)$.
When it is finite, the tilt $\mu^{(\lambda)}$ has density
$e^{\lambda y}/M_\mu(\lambda)$ relative to $\mu$.
Tonelli's theorem gives $M_{\mu*\mu'}(\lambda)=M_\mu(\lambda)M_{\mu'}(\lambda)$,
so tilting commutes with convolution. For an ID law $\mu$ with $M_\mu(\lambda)<\infty$, each
$n$th convolution root $\omega_n$ has a finite moment at $\lambda$, since $M_\mu(\lambda)=M_{\omega_n}(\lambda)^n$, and
$\mu^{(\lambda)}=(\omega_n^{(\lambda)})^{*n}$ for every $n$. Thus the tilt is ID
\cite[Example~33.14]{Sat99}.

Write $\beta=\re\rho$ and let $p_\theta$ be the density of the law
$\mu_\theta$ of $F_\theta$. Nakamura's density formula
\cite[Theorem~1.2]{Nak15a}, $p_\theta(y)=\theta([e^y]-e^y)e^{-\theta y}/\zeta(\theta)$, gives
\[
 e^{(\theta-\beta)y}p_\theta(y)
 =\frac{\theta\zeta(\beta)}{\beta\zeta(\theta)}p_\beta(y).
\]
In particular, $M_{\mu_\theta}(\lambda)=\theta\zeta(\beta)/(\beta\zeta(\theta))$ is finite, and the tilt of $\mu_\theta$ is $\mu_\beta$.
Thus the tilt of $\eta*\mu_\theta$ is $\eta^{(\lambda)}*\mu_\beta$.
Its characteristic function vanishes at $-{\im\rho}$
because $F_\beta(-{\im\rho})=0$. This rules out infinite
divisibility of the tilted law and hence of $\eta*\mu_\theta$, whose characteristic function is $HF_\theta$.
\end{proof}

By Hardy's theorem \cite{Har14} we may take $\rho$ on the critical line, so that $\lambda=\theta-\tfrac12$. Unlike Theorem~\ref{thm:main}, Theorem~\ref{thm:tilt} imposes no divisibility condition on $H$; its conclusion is weaker, since it excludes infinite divisibility but not quasi-infinite divisibility. Since a gamma law with rate $\theta$ has finite exponential moments at every $\lambda<\theta$, the choice $H=(\theta/(\theta-it))^{c-1}$ shows again that $\varphi_{\theta,c}$ is not ID. At $\theta=\tfrac12$ the moment condition is void, so $HF_{1/2}$ is not ID for any characteristic function $H$.

\section{A representation by prime powers}\label{sec:primes}

For $\theta>1$, the Euler product gives
\[
\log\frac{\zeta(\theta-it)}{\zeta(\theta)}
=\sum_p\sum_{k\ge1}\frac{p^{-k\theta}}{k}\bigl(e^{itk\log p}-1\bigr),
\]
a L\'evy--Khintchine formula whose L\'evy measure has an atom of mass
$p^{-k\theta}/k$ at each point $k\log p$. In the strip $0<\theta<1$ the Euler
product diverges, but a representation of the same shape survives once the
atoms are compensated by their average density. By the prime number theorem,
about $e^x\,\mathrm{d}x/x$ primes $p$ satisfy $\log p\in[x,x+\mathrm{d}x]$, and
each carries the weight $p^{-\theta}=e^{-\theta x}$, so the atoms have mean
density approximately $e^{(1-\theta)x}/x$ on $(0,\infty)$. The signed measure
$\kappa_\theta$ defined below consists of the atoms minus this mean density,
together with an elementary term on $(-\infty,0)$ that comes from the pole of
$\zeta$ at $s=1$. Its two parts on $(0,\infty)$ have infinite mass and cancel
only in an improper limit. The cancellation is controlled by the error term in
the prime number theorem, which is why a zero-free half-plane is needed. On each
interval $(\log m,\log(m+1))$ there are no prime powers, so only the negative
density remains; this is what excludes quasi-infinite divisibility.

Let $\Theta:=\sup\{\re\rho:\zeta(\rho)=0,\ 0<\re\rho<1\}$. The measure $\kappa_\theta$ below is defined for every $0<\theta<1$. In Theorem~\ref{thm:primes} we assume $\Theta<\theta<1$, so that $Z_\theta$ has no zeros, and write $\Psi_\theta$ for its continuous logarithm with $\Psi_\theta(0)=0$. The interval $(\Theta,1)$ is nonempty only if $\Theta<1$, which is not known; under the Riemann hypothesis $\Theta=\tfrac12$. Theorem~\ref{thm:primes} and Section~\ref{sec:numerics} are conditional in this sense, whereas Theorems~\ref{thm:main} and~\ref{thm:tilt}, Corollary~\ref{cor:xi}, and Proposition~\ref{prop:sharp} are unconditional. A proof of the Riemann hypothesis has recently been announced \cite{Alp26}; if it is confirmed, Theorem~\ref{thm:primes} and Section~\ref{sec:numerics} apply to every $\tfrac12<\theta<1$. Define a locally finite signed measure on the punctured real line by
\begin{equation}\label{eq:kernel}
\begin{split}
\kappa_\theta(\mathrm{d}x)
={}&\sum_p\sum_{k\ge1}\frac{p^{-k\theta}}{k}\,\delta_{k\log p}(\mathrm{d}x)\\
&-\one_{(0,\infty)}(x)\frac{e^{(1-\theta)x}}{x}\dd x
+\one_{(-\infty,0)}(x)\frac{e^{-(1-\theta)|x|}}{|x|}\dd x .
\end{split}
\end{equation}
Here $p$ runs over the primes. Both Jordan parts of $\kappa_\theta$ have infinite mass on $(1,\infty)$, by the divergence of $\sum_p p^{-1}$ and the exponential growth of the density.

We use the Fourier transform with the positive sign, $\F g(t)=\int_{\R}e^{itx}g(x)\dd x$ on the Schwartz space $\Sch(\R)$, and extend it by duality to tempered distributions.

\begin{theorem}\label{thm:primes}
Let $\Theta<\theta<1$. For every $t\in\R$,
\begin{equation}\label{eq:improper}
\begin{split}
\Psi_\theta(t)
={}&\int_{(-\infty,0)}(e^{itx}-1)\kappa_\theta(\mathrm{d}x)\\
&+\lim_{X\to\infty}\int_{(0,X]}(e^{itx}-1)\kappa_\theta(\mathrm{d}x).
\end{split}
\end{equation}
The limit uses a common cutoff for the atoms and the density on $(0,\infty)$. For every $h>0$, integration against $w_h\kappa_\theta$ defines a distribution on compactly supported smooth test functions on $\R$, which need not vanish near the origin. This distribution extends uniquely to a tempered distribution $W_{\theta,h}$. In $\Sch'(\R)$,
\begin{equation}\label{eq:Fourier-kernel}
D_h\Psi_\theta=\F W_{\theta,h}.
\end{equation}
\end{theorem}

\begin{proof}
The proof has four steps: an estimate for the error term $A$ in the prime number
theorem, a Laplace representation of $\log R$, the derivation of
\eqref{eq:improper}, and the passage to distributions. We use Frullani's formula
in the form
\begin{equation}\label{eq:Frullani}
\int_0^\infty\frac{e^{-ax}-e^{-bx}}{x}\dd x=\log\frac ba
\qquad(\re a>0,\ \re b>0),
\end{equation}
with the principal logarithm: for fixed $a$, both sides are holomorphic in $b$,
vanish at $b=a$, and have derivative $1/b$.

\emph{Step 1: the error term.} Put $\Pi(y):=\sum_{p^k\le y}1/k$, where $p$ runs
over the primes and $k\ge1$, and let $\pi(y)$ be the number of primes up to $y$.
For $x\ge0$ define
\[
E(x):=\int_0^x\frac{e^u-1}{u}\dd u,\qquad A(x):=\Pi(e^x)-E(x),
\]
and set $A(x)=0$ for $x<0$. For $x>0$ one has
$E(x)=\li(e^x)-\log x-\gamma_{\mathrm E}$, where $\li$ is the logarithmic
integral and $\gamma_{\mathrm E}$ is Euler's constant, so $A$ is the error term
in the prime number theorem for $\Pi$. Let $\psi(y):=\sum_{n\le y}\Lambda(n)$,
where $\Lambda(n)=\log p$ if $n=p^k$ and $\Lambda(n)=0$ otherwise. For every
$\varepsilon>0$, $\psi(y)=y+O_\varepsilon(y^{\Theta+\varepsilon})$
\cite[Lemma~4.2 and \S1]{Aka17}. Replacing $\psi(y)$ by $\sum_{p\le y}\log p$, or
$\Pi(y)$ by $\pi(y)$, changes it by $O(\sqrt y)$, because only prime powers with
$k\ge2$ are involved and Chebyshev's bounds apply; this is
$O(y^{\Theta+\varepsilon})$ because $\Theta\ge\tfrac12$ by Hardy's theorem
\cite{Har14}. Partial summation gives $\pi(y)-\li(y)=O_\varepsilon(y^{\Theta+\varepsilon})$,
hence $\Pi(y)-\li(y)=O_\varepsilon(y^{\Theta+\varepsilon})$, and therefore
\begin{equation}\label{eq:A-bound}
A(x)=O_\varepsilon(e^{(\Theta+\varepsilon)x})\qquad(x\to\infty).
\end{equation}

\emph{Step 2: a Laplace representation of $\log R$.} Define
$R(s):=(s-1)\zeta(s)/s$, with the removable value $R(1)=1$. The definition of
$\Theta$, the nonvanishing of $\zeta(1+it)$ for $t\ne0$ \cite[\S~3.3]{Tit86},
and the Euler product show that $R$ is holomorphic and zero-free for
$\re s>\Theta$. It is positive on $(\Theta,\infty)$, by the Euler product above
$1$ and the alternating Dirichlet series for $(1-2^{1-s})\zeta(s)$ below $1$. Let
$\log R$ be the holomorphic branch that is real on $(\Theta,\infty)$. For
$\re s>1$, the Euler product gives $\log\zeta(s)=\int_{[0,\infty)}e^{-sx}\dd\Pi(e^x)$,
and \eqref{eq:Frullani} gives
$\int_0^\infty e^{-sx}\dd E(x)=\int_0^\infty(e^{-(s-1)x}-e^{-sx})x^{-1}\dd x=\log\frac{s}{s-1}$.
Hence $\log R(s)=\int_{[0,\infty)}e^{-sx}\dd A(x)$ for $\re s>1$. Since $A(0)=0$,
integration by parts gives
\[
\int_{[0,X]}e^{-sx}\dd A(x)=e^{-sX}A(X)+s\int_0^X e^{-sx}A(x)\dd x .
\]
By \eqref{eq:A-bound}, the boundary term tends to zero for $\re s>\Theta$, and
the integral converges locally uniformly there. The identity theorem therefore
gives
\begin{equation}\label{eq:R}
\log R(s)=s\int_0^\infty e^{-sx}A(x)\dd x=\lim_{X\to\infty}\int_{[0,X]}e^{-sx}\dd A(x)
\qquad(\re s>\Theta).
\end{equation}

\emph{Step 3: proof of \eqref{eq:improper}.} For $s=\theta-it$,
\[
Z_\theta(t)=\frac{R(\theta-it)}{R(\theta)}\,
\frac{\theta-it}{\theta}\,\frac{1-\theta}{1-\theta+it}.
\]
The function $t\mapsto\log R(\theta-it)-\log R(\theta)$ and the principal
logarithms of the two elementary ratios are continuous and vanish at $t=0$, so
their sum is $\Psi_\theta$. By \eqref{eq:Frullani},
\begin{align*}
\log\frac{\theta-it}{\theta}
&=-\int_0^\infty(e^{itx}-1)\frac{e^{-\theta x}}{x}\dd x,\\
\log\frac{1-\theta}{1-\theta+it}
&=\int_{-\infty}^{0}(e^{itx}-1)\frac{e^{-(1-\theta)|x|}}{|x|}\dd x ,
\end{align*}
and by the second form of \eqref{eq:R},
$\log R(\theta-it)-\log R(\theta)=\lim_{X\to\infty}\int_{(0,X]}(e^{itx}-1)e^{-\theta x}\dd A(x)$.
On $(0,\infty)$, the measure $e^{-\theta x}\dd A(x)$ consists of the atoms
$p^{-k\theta}/k$ at $k\log p$ and the density $-(e^{(1-\theta)x}-e^{-\theta x})/x$.
Adding the density $-e^{-\theta x}/x$ of the first Frullani integral, which
converges absolutely and can therefore be moved inside the limit, gives the
restriction of $\kappa_\theta$ to $(0,\infty)$. The second Frullani integral is
the part of \eqref{eq:improper} on $(-\infty,0)$.

\emph{Step 4: the distributional identity.} Let $a_\theta(x):=e^{-\theta x}A(x)$.
By \eqref{eq:A-bound}, $a_\theta\in L^1(\R)$. The measure $e^{-\theta x}\dd A(x)$
is not tempered as it stands, because its density grows exponentially, but by
the product rule in $\mathcal{D}'(\R)$ it equals
$q_\theta:=\theta a_\theta+a_\theta'$, which is a tempered distribution. Since
$\F(a_\theta')(t)=-it\,\F a_\theta(t)$, the first form of \eqref{eq:R} gives
$\log R(\theta-it)=(\theta-it)\F a_\theta(t)=\F q_\theta(t)$; the constant
$\log R(\theta)$ is annihilated by $D_h$. For $g\in\Sch(\R)$, a direct
computation gives $\F(D_hg)=w_h\F g$. Since $D_h$ is symmetric under the dual
pairing and $w_h$ is smooth with bounded derivatives, it follows that
$D_h\F U=\F(w_hU)$ for $U\in\Sch'(\R)$, and hence
$D_h\log R(\theta-i\,\cdot)=\F(w_hq_\theta)$. On the two Frullani terms, $D_h$
acts under the integral sign, as in Lemma~\ref{lem:smoothing}, and multiplies
their densities by $w_h$; the products are integrable because $w_h(x)=O(x^2)$
at zero and the densities decay exponentially at infinity. On compactly
supported smooth test functions, $w_hq_\theta$ plus these two weighted densities
is exactly $w_h\kappa_\theta$. This proves \eqref{eq:Fourier-kernel}.
\end{proof}

Formula \eqref{eq:R} is a regularized version of the prime-power expansion in
Akatsuka \cite[Theorem~3]{Aka17}, specialized to the strict half-plane and to
$\log R$. For a related test-function identity under the Riemann hypothesis, see
Najnudel \cite[proof of Proposition~3.2]{Naj18}, where the pole term is given
explicitly.

The representation gives a second proof of Theorem~\ref{thm:main} for $\theta>\Theta$. If $Z_\theta$ were a quotient of ID characteristic functions, Lemma~\ref{lem:smoothing} would express $D_h\Psi_\theta$ as the Fourier transform of a finite signed measure. Injectivity on $\Sch'(\R)$ would identify that measure with $W_{\theta,h}$. On each interval
\[
I_m:=(\log m,\log(m+1)),\qquad m=1,2,\ldots,
\]
there are no prime-power atoms, since $p^k$ is an integer. The density there would be $-w_h(x)e^{(1-\theta)x}/x$. The intervals $I_m$ cover $(0,\infty)$ up to the countable set $\{\log m\}$, and $w_h(x)\to1$ as $x\to\infty$, so the total variation of $W_{\theta,h}$ would be at least $\int_1^\infty w_h(x)e^{(1-\theta)x}x^{-1}\dd x=\infty$.

For $\theta>1$, Nakamura \cite[Theorem~1.4 and its proof]{Nak15b} showed that the completed zeta distribution is QID but not ID: its quasi-L\'evy measure has a negative continuous density of finite variation away from zero. In the strip, the negative density has infinite variation at infinity, and this excludes quasi-infinite divisibility (Theorem~\ref{thm:main} and Corollary~\ref{cor:xi}).

\begin{remark}
For $\theta>1$ the same computation applies, but the elementary factor $(1-\theta)/(1-\theta+it)=(\theta-1)/(\theta-1-it)$ is the characteristic function of an exponential law on $(0,\infty)$. Frullani's formula therefore turns it into the density $e^{(1-\theta)x}/x$ on $(0,\infty)$ instead of a density on $(-\infty,0)$. This cancels the negative density and leaves only the prime-power atoms, as in the Euler product at the beginning of this section.
\end{remark}

\begin{figure}[H]
\centering
\includegraphics[width=\linewidth]{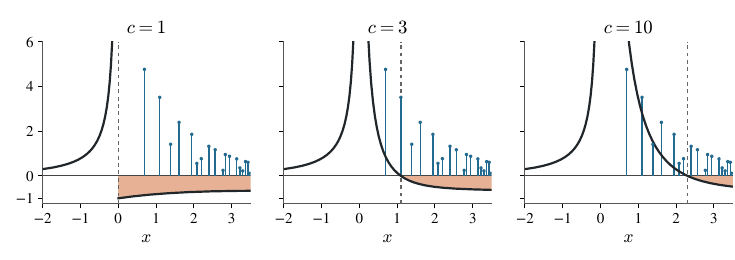}
\caption{The signed measure $\kappa_{\theta,c}$ for $\theta=0.75$ and $c=1,3,10$. The curve is its continuous density; the stems at $k\log p$ are its atoms, drawn with height $8p^{-k\theta}/k$ (scaled by $8$ for visibility). The shaded region marks where the density is negative, and the dashed lines mark $x=\log c$.}
\label{fig:signed-kernel}
\end{figure}

For $\Theta<\theta<1$ and $c\ge1$, the signed kernel for the distinguished logarithm of $\varphi_{\theta,c}$ is given by \eqref{eq:Frullani} as
\[
\kappa_{\theta,c}(\mathrm{d}x)
:=\kappa_\theta(\mathrm{d}x)
+c\,\one_{(0,\infty)}(x)\frac{e^{-\theta x}}{x}\dd x .
\]
The density of its continuous part on $(0,\infty)$ is
$(c\,e^{-\theta x}-e^{(1-\theta)x})/x$. It is negative for $x>\log c$, hence throughout $I_m$ when $m\ge c$. Thus the gamma factor compensates the negative density only on $(0,\log c)$, and the argument following the proof of Theorem~\ref{thm:primes} applies verbatim to $\varphi_{\theta,c}$. Figure~\ref{fig:signed-kernel} shows the continuous density and the prime-power atoms.

The hypothesis $\theta>\Theta$ in Theorem~\ref{thm:primes} is sharp, except possibly at $\theta=\Theta$. The measure $\kappa_\theta$ is defined for every $0<\theta<1$, and the functions $A$ and $R$ in the proof of Theorem~\ref{thm:primes} depend neither on $\theta$ nor on $\Theta$.

\begin{proposition}\label{prop:sharp}
Let $0<\theta<1$, and suppose that the limit in \eqref{eq:improper} exists for some $t\ne0$. Then $\theta\ge\tfrac12$, $\zeta(s)\ne0$ for $\re s>\theta$, and the exponential of the right-hand side of \eqref{eq:improper} equals $Z_\theta(t)$. Consequently, the limit in \eqref{eq:improper} exists for every $t\in\R$ if $\Theta<\theta<1$, and for no $t\ne0$ if $0<\theta<\Theta$.
\end{proposition}

\begin{proof}
On $(0,\infty)$ we have $\kappa_\theta(\mathrm{d}x)=e^{-\theta x}\dd A(x)-e^{-\theta x}\dd x/x$, and $(e^{itx}-1)e^{-\theta x}/x$ is integrable there. Hence
\[
B(X):=\int_{(0,X]}(e^{itx}-1)e^{-\theta x}\dd A(x)
\]
has a finite limit $B(\infty)$ as $X\to\infty$, and $B$ is bounded on $[0,\infty)$. Integration by parts shows that, locally uniformly for $\re s>\theta$,
\[
\lim_{X\to\infty}\int_{(0,X]}(e^{itx}-1)e^{-sx}\dd A(x)=\Phi(s):=(s-\theta)\int_0^\infty e^{-(s-\theta)x}B(x)\dd x ,
\]
so $\Phi$ is holomorphic for $\re s>\theta$. For $\re s>1$, as in Step~2 of the proof of Theorem~\ref{thm:primes}, $e^{\Phi(s)}=R(s-it)/R(s)$. The function $R$ is holomorphic for $\re s>0$ and positive on $(0,\infty)$, so the identity theorem shows that $R(s-it)/R(s)$ is holomorphic and zero-free for $\re s>\theta$. Thus the zeros of $\zeta$ in $\re s>\theta$, counted with multiplicity, are invariant under $s\mapsto s+it$. Letting $s\downarrow\theta$ gives $\Phi(s)\to B(\infty)$, hence $R(\theta-it)=R(\theta)e^{B(\infty)}\ne0$. The right-hand side of \eqref{eq:improper} equals $B(\infty)+\log((\theta-it)/\theta)+\log((1-\theta)/(1-\theta+it))$ by Step~3 of the proof of Theorem~\ref{thm:primes}, and its exponential is $Z_\theta(t)$.

Let $N(V)$ be the number of zeros $\rho$ of $\zeta$ with $0<\re\rho<1$ and $0<\im\rho\le V$, and let $N(\sigma,V)$ be the number of those with $\re\rho>\sigma$, all counted with multiplicity. We use $N(V)\sim(V/2\pi)\log V$ \cite[Theorem~9.4]{Tit86} and, for fixed $\tfrac12<\sigma<1$, $N(\sigma,V)=o(V)$ \cite[Theorem~9.17]{Tit86}. Suppose first that $\theta<\tfrac12$. The rectangle $\theta\le\re s\le1$, $0\le\im s\le|t|$ contains finitely many zeros, so translation invariance shows that at most $O(V)$ zeros satisfy $\re\rho>\theta$ and $0<\im\rho\le V$. By the functional equation \cite[Chap.~II]{Tit86}, the map $\rho\mapsto1-\bar\rho$ sends the remaining zeros to zeros with real part at least $1-\theta$ and the same imaginary part. For $\tfrac12<\sigma<1-\theta$ their number is therefore at most $N(\sigma,V)=o(V)$, which contradicts $N(V)\sim(V/2\pi)\log V$. Hence $\theta\ge\tfrac12$. Suppose next that $\zeta(\rho)=0$ with $\re\rho>\theta$. Then $\rho+int$ is a zero for every $n\in\mathbb{Z}$, so $N(\sigma,V)\ge V/|t|-1$ for $\tfrac12<\sigma<\re\rho$. This again contradicts $N(\sigma,V)=o(V)$.

The final assertion follows from Theorem~\ref{thm:primes} for $\Theta<\theta<1$, and from what was just proved for $0<\theta<\Theta$.
\end{proof}

\begin{remark}
For the integrals without differences, $\int_{(0,X]}e^{-sx}\dd A(x)$, convergence at one point of a line $\re s=\sigma_0$ with $\tfrac12<\sigma_0<1$ already forces $\zeta(s)\ne0$ for $\re s>\sigma_0$, by Akatsuka \cite[Theorem~2 and Lemma~2.2]{Aka17}; the Abelian step in the proof above follows \cite[Lemma~3.3]{Aka17}. Proposition~\ref{prop:sharp} assumes only the convergence of the differences $\int_{(0,X]}(e^{itx}-1)e^{-\theta x}\dd A(x)$, which is what \eqref{eq:improper} involves, and it also covers $0<\theta\le\tfrac12$.
\end{remark}

\section{Numerical illustration}\label{sec:numerics}

Theorem~\ref{thm:primes} and the discussion of $\kappa_{\theta,c}$ that follows
its proof predict that the signed measure $\kappa_{\theta,c}$, smoothed and
computed from values of $\zeta$ on the line $\re s=\theta$ alone, takes negative
values, which an ID characteristic function cannot produce. This section checks the representation and the negativity
numerically. It is an illustration, not part of any proof, and it assumes
$\theta>\Theta$; for the values of $\theta$ used below this requires $\Theta<0.6$.

\emph{The test.} Let $h,T>0$, and let $g_T$ be the centered Gaussian density with
standard deviation $1/T$, whose characteristic function is $e^{-t^2/(2T^2)}$. For
a zero-free characteristic function $\chi$ with distinguished logarithm
$\log\chi$, define, whenever the integral converges absolutely,
\begin{equation}\label{eq:K}
K_{T,h}[\chi](x):=\frac{1}{2\pi}\int_\R D_h\log\chi(t)\,e^{-itx-t^2/(2T^2)}\dd t .
\end{equation}
This is the inverse Fourier transform of $D_h\log\chi$, smoothed at scale $1/T$.
If $\chi$ is ID with Gaussian coefficient $a$ and L\'evy measure $\nu$,
Lemma~\ref{lem:smoothing} gives
$D_h\log\chi=\F\bigl(\tfrac{ah^2}{6}\delta_0+w_h\nu\bigr)$, and hence
$K_{T,h}[\chi]=g_T*\bigl(\tfrac{ah^2}{6}\delta_0+w_h\nu\bigr)\ge0$. Negative values
of $K_{T,h}[\chi]$ therefore show that $\chi$ is not ID.

\emph{The prediction.} Let $\chi=\varphi_{\theta,c}$ with $\Theta<\theta<1$ and
$c\ge1$. By \eqref{eq:R} and the elementary factors, $\log\chi(t)=O(1+|t|)$, so
\eqref{eq:K} converges absolutely, and Theorem~\ref{thm:primes}, together with
\eqref{eq:Frullani} for the gamma factor, gives
\begin{equation}\label{eq:K-pred}
K_{T,h}[\varphi_{\theta,c}]=g_T*(w_h\kappa_{\theta,c}).
\end{equation}
The convolution is an ordinary integral, because the Gaussian decays faster
than the variation of $w_h\kappa_{\theta,c}$ on unit intervals grows, which is
at most exponentially. As $T\to\infty$, the right-hand side converges, uniformly
on compact subintervals of each gap $I_m=(\log m,\log(m+1))$, an interval that contains no prime powers, to the density
$-w_h(x)\bigl(e^{(1-\theta)x}-c\,e^{-\theta x}\bigr)/x$, which is negative when
$m\ge c$. Indeed, splitting $w_h\kappa_{\theta,c}$ with a smooth cutoff supported
in $I_m$, the part inside $I_m$ is a continuous density, on which $g_T$ acts as
an approximate identity, and the convolution of $g_T$ with the remaining part
tends to zero on the compact subinterval by Gaussian decay.

\begin{table}[!htb]
\centering
\begin{tabular}{ccrccc}
\hline
$\theta$ & $c$ & Minimum & $x_{\min}$ & Gap & Relative error\\
\hline
$0.60$ & $1$  & $-1.257$ & $3.804$ & $I_{44}$ & $1.2\times10^{-5}$\\
$0.60$ & $3$  & $-1.199$ & $3.804$ & $I_{44}$ & $1.2\times10^{-5}$\\
$0.60$ & $10$ & $-1.008$ & $3.928$ & $I_{50}$ & $1.1\times10^{-5}$\\
$0.75$ & $1$  & $-0.922$ & $1.444$ & $I_{4}$  & $1.1\times10^{-5}$\\
$0.75$ & $3$  & $-0.686$ & $3.664$ & $I_{39}$ & $1.1\times10^{-5}$\\
$0.75$ & $10$ & $-0.563$ & $3.804$ & $I_{44}$ & $1.1\times10^{-5}$\\
$0.90$ & $1$  & $-0.750$ & $1.324$ & $I_{3}$  & $1.0\times10^{-5}$\\
$0.90$ & $3$  & $-0.400$ & $3.524$ & $I_{33}$ & $1.0\times10^{-5}$\\
$0.90$ & $10$ & $-0.319$ & $3.664$ & $I_{39}$ & $1.1\times10^{-5}$\\
\hline
\end{tabular}
\caption{Minimum of $\re K^{\mathrm{num}}$ over $\mathcal X$, the point
$x_{\min}$ where it is attained and the gap $I_m=(\log m,\log(m+1))$ containing
it, and the relative error with respect to the prime-power prediction
\eqref{eq:K-pred}; $T=80$, $h=3$.}
\label{tab:numerics}
\end{table}

\emph{The computation.} For $\theta\in\{0.6,0.75,0.9\}$, $c\in\{1,3,10\}$,
$h=3$, and $T=80$, we computed the two sides of \eqref{eq:K-pred} independently,
on the grid $\mathcal X$ from $-2$ to $4$ with step $0.004$.
\begin{itemize}
\item The left-hand side, $K^{\mathrm{num}}$, uses only values of $\zeta$. We
sampled $\zeta(\theta-it)$ in double precision with \texttt{mpmath} for
$0\le t\le480$ at step $0.02$, extended to $t<0$ by conjugation, obtained
$\log Z_\theta$ by unwrapping the phase from $t=0$, and subtracted
$c\log((\theta-it)/\theta)$. We computed $D_3$ by the trapezoidal rule, which
leaves $|t|\le477$, and \eqref{eq:K} by a Riemann sum; the Gaussian factor at
$|t|=477$ is below $e^{-17}$.
\item The right-hand side, $K^{\mathrm{pred}}$, uses only the primes: the atoms
of $\kappa_{\theta,c}$ at prime powers $p^k\le64$, and its continuous density by
quadrature on $[-6,7)$ with step $0.0002$.
\end{itemize}
The code is available on request.

\emph{Results.} Table~\ref{tab:numerics} reports, for each case, the minimum of
$\re K^{\mathrm{num}}$ over $\mathcal X$, the point where it is attained, and
the relative error
\[
\max_{\mathcal X}\bigl|\re K^{\mathrm{num}}-K^{\mathrm{pred}}\bigr|\Big/\max_{\mathcal X}\bigl|K^{\mathrm{pred}}\bigr| .
\]
By symmetry $K^{\mathrm{num}}$ is real; its computed imaginary part is below
$10^{-11}$ in all cases. The two computations agree to about five digits in all
nine cases, which supports \eqref{eq:improper}. All minima are negative, and each is attained inside
a gap $I_m$ with $m\ge c$, as \eqref{eq:K-pred} predicts. Most of them lie near
the right end of the grid, because the negative density grows like
$e^{(1-\theta)x}/x$. For an ID law the minimum would be nonnegative, so the
table illustrates that $\varphi_{\theta,c}$ is not ID, even for $c=10$.

\subsection*{Use of AI tools}
During the preparation of this work the author used Claude (Anthropic) and ChatGPT/Codex (OpenAI) to draft and revise the text and to assist with the numerical computations and the figure. After using these tools, the author reviewed and edited the content as needed and takes full responsibility for the content of the published article.


\begin{thebibliography}{99}

\bibitem{Aka17}
H. Akatsuka, \emph{The Euler product for the Riemann zeta-function in the critical strip},
Kodai Math. J. \textbf{40} (2017), no.~1, 79--101.

\bibitem{Alp26}
L. Alp\"oge, \emph{Remark}, announcement of a proof of the Riemann hypothesis with a Lean formalization,
November 20, 2026. Available at \texttt{https://alpo.ge/remark.pdf}.

\bibitem{FM26}
A. Feuerverger and G. Martin, \emph{Randomness and the Riemann hypothesis: Probability, statistics, and the primes},
Miscellaneous Books, vol.~156, American Mathematical Society, Providence, RI, 2026.

\bibitem{Har14}
G. H. Hardy, \emph{Sur les z\'eros de la fonction $\zeta(s)$ de Riemann},
C. R. Acad. Sci. Paris \textbf{158} (1914), 1012--1014.

\bibitem{LPS18}
A. Lindner, L. Pan, and K. Sato, \emph{On quasi-infinitely divisible distributions},
Trans. Amer. Math. Soc. \textbf{370} (2018), no.~12, 8483--8520.

\bibitem{LS52}
E. Lukacs and O. Sz\'asz, \emph{On analytic characteristic functions},
Pacific J. Math. \textbf{2} (1952), no.~4, 615--625.

\bibitem{Naj18}
J. Najnudel, \emph{On the extreme values of the Riemann zeta function on random intervals of the critical line},
Probab. Theory Related Fields \textbf{172} (2018), nos.~1--2, 387--452.

\bibitem{Nak15a}
T. Nakamura, \emph{A modified Riemann zeta distribution in the critical strip},
Proc. Amer. Math. Soc. \textbf{143} (2015), no.~2, 897--905.

\bibitem{Nak15b}
T. Nakamura, \emph{A complete Riemann zeta distribution and the Riemann hypothesis},
Bernoulli \textbf{21} (2015), no.~1, 604--617.

\bibitem{Nak16}
T. Nakamura, \emph{Hurwitz--Lerch zeta and Hurwitz--Lerch type of Euler--Zagier double zeta distributions},
Infin. Dimens. Anal. Quantum Probab. Relat. Top. \textbf{19} (2016), no.~4, 1650029, 12 pp.

\bibitem{Sat99}
K. Sato, \emph{L\'evy processes and infinitely divisible distributions},
Cambridge Studies in Advanced Mathematics, vol.~68, Cambridge University Press, Cambridge, 1999.

\bibitem{Tit86}
E. C. Titchmarsh, \emph{The theory of the Riemann zeta-function},
2nd ed., revised by D. R. Heath-Brown, Clarendon Press, Oxford, 1986.

\bibitem{Vor75}
S. M. Voronin, \emph{Theorem on the \textquotedblleft universality\textquotedblright\ of the Riemann zeta-function},
Math. USSR-Izv. \textbf{9} (1975), no.~3, 443--453.

\end{thebibliography}
\end{document}